\input amssym.tex

\mag=1200 
\hsize=130mm  \vsize=180mm  \voffset 5mm
\lineskiplimit=8pt \lineskip=8pt
\tolerance=10000 \pretolerance=1000 \parindent=0mm \raggedright

 \font\uq =cmr10 at  14pt 
 \font\uc =cmr10 at  15pt

\font\tengo=eufm10

\font\sevengo=eufm7

\font\fivego=eufm5

\font\tenbb=msbm7 at 10pt

\font\sevenbb=msbm7   
 
\font\fivebb=msbm5


\newfam\gofam  \textfont\gofam=\tengo
\scriptfont\gofam=\sevengo   \scriptscriptfont\gofam=\fivego
\def\go{\fam\gofam\tengo}

\newfam\bbfam  \textfont\bbfam=\tenbb
\scriptfont\bbfam=\sevenbb   \scriptscriptfont\bbfam=\fivebb
\def\bb{\fam\bbfam\tenbb}


\def\ind{\hskip 1em\relax}

\let\mbox=\hbox

\def\dim{\mathop{\rm dim}\nolimits}

\def\Ext{\mathop{\rm Ext}\nolimits}
\def\R{{\bb R}} \def\C{{\bb C}}  \def\Z{{\bb Z}}

\def\cl{\centerline}
\def\{{\lbrace}  
\def\}{\rbrace}   
\def\({\langle}  
\def\){\rangle}
\def\[{\lbrack} 
\def\]{\rbrack}
\def\l{\ell}
\def\arrow{\rightarrow}
\def\iso{\buildrel\sim\over{\arrow}}
\def\.{\bullet}
\def\bs{\bigskip}

\def\Spin{\hbox{Spin}}

\def\n{\cap}

\def\\{\backslash}
\def\back/{\backslash}

%
%
\cl{\uc A simple question about a complicated object} \bs
%
%
%
The complicated object is the cohomological induction functor (for which the 
biblical reference is Knapp \& Vogan [KV]). Here is the simple question. For 
definiteness sake set 
$G := PU(n,1)$, \hbox{$H := PU(k,1)$} with $k \leq n$, and let 
${\cal H}_\rho^n$ (resp. 
${\cal H}_\rho^k$) be the category of Harish-Chandra modules with the 
generalized 
infinitesimal character of the trivial module attached to $G$ (resp. $H$). 
I'll surprise nobody by claiming that $H$ is a subgroup of $G$. What's far 
less obvious, 
but proved by Khoroshkin [K], is the existence of a full embedding $F$ of 
${\cal H}_\rho^k$ into ${\cal H}_\rho^n$ (Fuser [F,Thm I.4.2] showed that $F$ 
is even
Ext-full, that is compatible with Ext calculus), prompting the question : is 
there a geometric interpretation of the embedding $F$ ? The first candidate 
for $F$ is the 
(ordinary) induction functor~; but this fails miserably --- 
so let's break the pseudo-suspense of this introduction by saying that I 
claim that 
$F$ is (isomorphic to) a certain cohomological induction functor, and 
conjecture that this phenomenon is general. \bs
\cl{\uq 1. \ Statements} \bs
Let $G$ be a center free connected semisimple Lie group, $K \subset G$ a 
maximal 
compact subgroup, $\go g \supset \go k$ the respective complexified Lie 
algebras. 
Let's start by recalling the notion of Harish-Chandra module. 
Say that a $\go g$-module $V$ is $\go k$-{\bf finite} if it is a sum of 
finite dimensional sub-$\go k$-modules, and that $V$ is an 
\hbox{$(\go g,k)$-{\bf module}} if 
it is $\go k$-finite and \hbox{$\go k$-semisimple.} The category 
${\cal H = H}({\go g},K)$ of {\bf Harish-Chandra modules} is the full subcategory 
of 
\hbox{$\go g$-mod} whose objects are those $(\go g,k)$-modules of finite 
length $V$ 
such that for any finite dimensional \hbox{$\go k$-invariant} subspace 
$F \subset V$ the action of $\go k$ on $F$ exponentiates to $K$. 
The category $\cal H$ is a $\C$-category in the sense of Bass [B] page~57. 
Let $I$ be the annihilator of the trivial module in the center of $U(\go g)$, 
let 
$${\cal H_\rho = H}_\rho({\go g},K)$$
be the full sub-$\C$-category of $\cal H$ whose objects are annihilated by 
some power of $I$, let $\cal I$ be 
the set of isomorphism classes of simple objects of $\cal H_\rho$ [it is a 
finite 
set]~; for each $i \in \cal I$ choose a representative $V_i \in i$ and let
$\l(i)$ be the {\bf projective dimension} of $V_i \ $ [{\it i.e.} the supremum in 
$\Z \cup \{+ \infty\}$ of the set 
$\{n \in \Z \ | \ \Ext^n(V_i \, ,-) \not= 0 \}$]. \bs 
%
%
\vbox{ {\parindent=6mm 
\item{(1)} {\bf Definition}. The {\bf $\cal H_\rho$-ordering} is the smallest 
partial ordering $\leq$ on $\cal I$ satisfying \bs 
$$\left.\matrix{i, j \in \cal I \ \cr \cr
\l(j) = \l(i) + 1 < \infty \ \cr \cr
\Ext^1(V_j,V_i) \not= 0 \ }
\right\} \ \Longrightarrow \ i \leq j.$$} \par} \bs  
%
%
{\parindent=6mm 
\item{(2)} {\bf Definition}. The {\bf sub-$\C$-category generated by the subset} 
$\cal J$ {\bf of} $\cal H_\rho$ is the full sub-$\C$-category 
$\( \cal J \)_{H_\rho}$ of 
$\cal H_\rho$ characterized by the condition that an object $V$ of 
$\cal H_\rho$ 
belongs to $\( \cal J \)_{H_\rho}$ iff each simple subquotient of $V$ is 
isomorphic to 
$V_j$ for some $j \in \cal J$. \par}   \bs 
%
%
%
{\parindent=6mm 
\item{(3)}  {\bf Definition}. Say that a full sub-$\C$-category $\cal C$ of 
$\cal H_\rho$ is {\bf Ext-full} in $\cal H_\rho$ if for all $V, W \in \cal C$ 
the natural morphism 
$$\Ext_{\cal C}^\.(V,W) \rightarrow \Ext_{{\cal H}_\rho}^\.(V,W)$$
is an isomorphism. \par} \bs 
%
%
For $i \in \cal I$ put ${\cal J}_i := \{ j \in {\cal I} \ | \ j \leq i \}$, 
let $\theta$ be the Cartan involution of $({\go g},K)$, denote by 
$d$ the dimension of $G/K$, and consider the following  \bs
%
%
\vbox{{\parindent=6mm
\item{(4)} {\bf Property of} $G$. For each \hbox{$i \in \cal I$} such that 
$V_i$ is unitary the cohomology $H^{d- \l(i)}({\go g}, K \, ; V_i)$ is 
nonzero and there is a $\theta$-stable parabolic subalgebra of $\go g$ with 
Levi subgroup $L = L_i$ (see Vogan [V2,4.1,4.2] for definitions) 
satisfying \medskip
\itemitem{(a)} the corresponding cohomological induction functor $F$ 
(see [KV]) sets up an equivalence
$${\cal H}_\rho({\go l},L \n K) \ \iso \ 
 \({\cal J}_i \)_{ {\cal H}_\rho({\go g}, K)} \ ; $$
\itemitem{(b)} $F \, \C \simeq V_i \ ;$  \medskip
\itemitem{(c)} $\({\cal J}_i \)_{ {\cal H}_\rho({\go g}, K)}$ is Ext-full 
in ${\cal H}_\rho({\go g},K)$~; \medskip
\itemitem{(d)} if $a$ is nonzero vector of $H^{d- \l(i)}({\go g}, K \, ; V_i)$
and $V$ a simple object of ${\cal H}_\rho({\go l},L \n K)$, then the map 
$$\matrix{
H^\.({\go l},L \n K ; V) & \arrow & H^{d- \l(i)+ \.}({\go g},K ; FV) \cr \cr
x & \mapsto & F(x) \cup a  \cr}$$
[where $\cup$ denotes the cup-product] is an isomorphism [of 
$H^\.({\go l},L \n K ;\C)$-modules]~; \bs
\itemitem{(e)} we have $2 \, \l(i) \ = \ d \ + \ \dim \ L/(L \n K)$. \par}} \bs
Note once and for all that (e) follows from (a) by the well known argument 
which consists in setting $V := \C$ and using Poincar\'e duality. \bs
{\parindent=6mm
\item{(6)} {\bf Conjecture}. 
All center free connected semisimple Lie groups have Property~(4). 
\par}  \bs
\vbox{\ind A partial proof (with explicitly indicated gaps) of the fact that 
$PU(n,1)$, $P \Spin(n,1)$ and $SL(3, \R)$ have Property~(4) is contained 
in the expanded version of this text, downloadable from \bs
\cl{http://www.iecn.u-nancy.fr/$\sim$gaillard/Recherche/Ci/ci.html}} \bs
%
%
\cl{* \ * \ *} \bs
{\parindent=10mm
\item{[F]} Fuser A., Autour de la conjecture d'Alexandru, 
Th\`ese de l'Univesit\'e Nancy 1 (1997). \bs
\item{[K]} Khoroshkin S.M., Category of Harish-Chandra modules of the group 
SU(n,1).
{\it Funct. Anal. Appl.} {\bf 14} (1980) 153-155. \bs
\item{[KV]} Knapp A. \& Vogan D., {\bf Cohomological induction
and unitary representations}, Princeton University Press (1995). \bs
\item{[V]} Vogan D., Cohomology and group representations, 
{\it Proc. Symp. Pure Math.} {\bf 61}  (1997) 219-243. \bs
\par} 
\bye